\documentclass{article}

\usepackage[utf8]{inputenc}
\usepackage[T1]{fontenc}
\usepackage{graphicx}
\usepackage{tikz}
\usepackage{amsmath,mathtools,amsfonts,amsthm}
\usepackage[pdftex,colorlinks=true,linkcolor=blue,citecolor=blue,urlcolor=blue]{hyperref}
\usepackage{csquotes}
\usepackage{tikz}
\usepackage{bm}

\newcommand\y{\bm{y}}
\newcommand\lam{\bm{\lambda}}
\newcommand\R{\mathbb{R}}
\renewcommand\d{\mathrm{d}}
\renewcommand\L{\mathcal{L}}

\newtheorem{remark}{Remark}
\newtheorem{theorem}{Theorem}

\title{Non-overlapping Schwarz methods in time for parabolic optimal control problems}

\author{Martin Jakob Gander$^{1}$, Liu-Di LU$^{1}$}
\date
{%
\noindent{\small\textit{$^1$Section of Mathematics, University of Geneva, Rue du Conseil Général 7-9, 1205 Geneva, Switzerland}}\\
}

\begin{document}

\maketitle

\begin{abstract} 
We present here the classical Schwarz method with a time domain decomposition applied to unconstrained parabolic optimal control problems. Unlike Dirichlet--Neumann and Neumann--Neumann algorithms, we find different properties based on the forward-backward structure of the optimality system. Variants can be found using only Dirichlet and Neumann transmission conditions. Some of these variants are only good smoothers, while others could lead to efficient solvers.
\end{abstract}

\section{Introduction}\label{sec:1}

The classical Schwarz method, originally introduced by Hermann Amandus
Schwarz to prove existence and uniqueness of solutions to Laplace's
equation~\cite{Schwarz1870}, has since then been extensively studied
and applied to a wide range of problems. A historical review can be
found in~\cite{Gander2008}. It is well known that the method fails to
converge when applied to non-overlapping spatial subdomains due to the
repeated passing of identical "Dirichlet data" from one subdomain to
the other. Several modified methods have then been proposed to address
this issue, notably by Lions~\cite{Lions1990}. More recently, Schwarz
methods in time have been proposed for the time-parallel solution of
parabolic optimal control problems, as discussed
in~\cite{Gander2016}, and it was noted:
\begin{displayquote}
  We present a rigorous convergence analysis for the case of two
  subdomains, which shows that the classical Schwarz method converges,
  even without overlap! Reformulating the algorithm reveals that this
  is because imposing initial conditions for $y$ and final conditions
  on $\lambda$ is equivalent to using Robin transmission conditions
  between time subdomains for y.
\end{displayquote} 
To gain deeper insight into the convergence of the classical Schwarz
method applied to parabolic optimal control problems with a
non-overlapping time domain decomposition, we study the following
model problem: for a given desired state $\hat \y(t)$ and parameters
$\gamma,\nu > 0$, we want to solve the minimization problem
\begin{equation}\label{eq:model}
\begin{aligned}
\min_{\y, u}&\frac12 \int_0^T \|\y -\hat \y\|^2 \d t  + \frac{\gamma}2 \|\y(T) - \hat \y(T)\|^2 + \frac{\nu}2 \int_0^T \|u\|^2 \d t, \\
&\text{ subject to } \quad \dot \y + A \y = \mathbf{u}, \quad \y(0) = \y_0, 
\end{aligned}
\end{equation}
where $\dot \y + A \y = \mathbf{u}$ represents the semi-discretization
of a parabolic partial differential equation (PDE) of the form
$\partial_t y + \L y = u$. Here, $\mathbf{u}$ is the control variable,
and $\y_0$ denotes the initial condition. By applying the Lagrange
multiplier approach and eliminating the control variable $\mathbf{u}$,
we derive from~\eqref{eq:model} the first-order optimality system
(see, e.g.,~\cite[Section 2]{Gander20241})
\begin{equation}\label{eq:forbacksys}
\left\{
\begin{aligned}
\begin{pmatrix}
\dot \y \\
\dot \lam
\end{pmatrix}
+ 
\begin{pmatrix}
A & -\nu^{-1} I\\
-I & -A^T
\end{pmatrix}
\begin{pmatrix}
\y \\
\lam
\end{pmatrix}
&=
\begin{pmatrix}
0\\
-\hat \y
\end{pmatrix}
\quad \text{ in } (0, T), \\
\y(0) &= \y_0, \\
\lam(T) + \gamma \y(T) &= \gamma \hat \y(T),
\end{aligned}
\right.
\end{equation}
where $\lam$ is the adjoint state. The system described
by~\eqref{eq:forbacksys} is a forward-backward system of ordinary
differential equations (ODEs), in which the state $\y$ propagates
forward in time starting from an initial condition, while the adjoint
state $\lam$ propagates backward in time with a final condition.

To investigate the application of a non-overlapping classical Schwarz
method in time to solve this system, we decompose the time interval
$(0, T)$ into two non-overlapping subdomains: $I_1 := (0, \alpha)$ and
$I_2 := (\alpha, T)$, where $\alpha$ represents the interface. In
Section~\ref{sec:2}, we first introduce four variants of the classical
Schwarz algorithm and analyze their convergence behavior. In
Section~\ref{sec:3}, we replace the Dirichlet transmission condition
used in the classical Schwarz algorithm with a Neumann transmission
condition and study the resulting convergence properties. Finally, we
discuss our results in Section~\ref{sec:4} and conclude with some
comments.

\section{Dirichlet transmission conditions}\label{sec:2}

\begin{figure}[t]
  \centering
  \begin{tikzpicture}
  \draw[thick, ->] (0, 0)--(4.3, 0) node[anchor = north west] {$x$};
  \draw[thick, ->] (0, 0)--(0, 2.3) node[anchor = south east] {$t$};
  \draw[dashed] (0, 2)--(4, 2);
  \draw[dashed] (4, 2)--(4, 0);
  \draw[dashed] (2, 0)--(2, 2);
  \node at (4, -0.15) {$L$};
  \node at (0.05, -0.15) {0};
  \node at (-0.15, 0.05) {0};
  \node at (-0.15, 2.00) {$T$};
  \node at (2, -0.15) {$\alpha$};
  \draw[thick, dashed, <->] (1.7, 1.5)--(2.3, 1.5);
  \draw[thick, dashed, <->] (1.7, 0.5)--(2.3, 0.5);
  \node at (1.9, 0.3) {$y$};
  \node at (2.1, 1.7) {$\lambda$};
  \end{tikzpicture}
  \begin{tikzpicture}
  \draw[thick, ->] (0, 0)--(4.3, 0) node[anchor = north west] {$x$};
  \draw[thick, ->] (0, 0)--(0, 2.3) node[anchor = south east] {$t$};
  \draw[dashed] (0, 2)--(4, 2);
  \draw[dashed] (4, 2)--(4, 0);
  \node at (4, -0.15) {$L$};
  \node at (0.05, -0.15) {0};
  \node at (-0.15, 0.05) {0};
  \node at (-0.15, 2.00) {$T$};
  \draw[thick, dashed, ->] (1, 0)--(1, 0.6);
  \draw[thick, dashed, ->] (3, 2)--(3, 1.4);
  \node at (0.7, 0.3) {$y$};
  \node at (3.3, 1.7) {$\lambda$};
  \draw[dashed] (0, 1)--(4, 1);
  \node at (-0.15, 1.00) {$\alpha$};
  \draw[thick, dashed, ->] (1, 1)--(1, 1.6);
  \draw[thick, dashed, ->] (3, 1)--(3, 0.4);
  \node at (0.7, 1.3) {?};
	\node at (3.3, 0.7) {?};
  \end{tikzpicture}
  \caption{One dimensional illustration of decomposing in space (left)
    and decomposing in time (right).}
  \label{fig:illustration}
\end{figure}
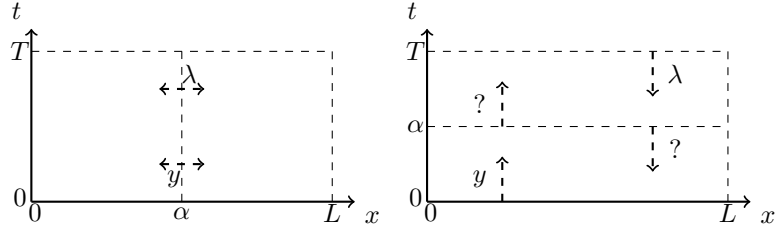

When decomposing in space, the standard way is to pass the values for
the state $y$ and the adjoint state $\lambda$ from one subdomain to
its neighbor, as illustrated in Figure~\ref{fig:illustration} on the
left.
However, this becomes much more tricky when decomposing in time as
shown in Figure~\ref{fig:illustration} on the right. Since the
system~\eqref{eq:forbacksys} is a forward-backward system, it
initially seems natural to preserve this property in the decomposed
case, and to transmit in $I_1$ a final condition for the adjoint state
$\lambda$, while an initial condition for $y$ is already
present. Similarly, in $I_2$, where a final condition for $\lambda$ is
already present, it is natural to transmit an initial condition for
the state $y$. A natural Schwarz algorithm in time hence solves for
the iteration index $k = 1, 2, \ldots$
\begin{equation}\label{eq:SD1}
\begin{aligned}
\left\{
\begin{aligned}
\begin{pmatrix}
\dot \y^k_1\\
\dot \lam^k_1
\end{pmatrix}
+ 
\begin{pmatrix}
A & -\nu^{-1} I\\
-I & -A^T
\end{pmatrix}
\begin{pmatrix}
\y^k_1\\
\lam^k_1
\end{pmatrix}
&=
\begin{pmatrix}
0\\
-\hat \y
\end{pmatrix}
\quad \text{ in } (0, \alpha), \\
\y^k_1(0) &= \y_0, \\
\lam^k_1(\alpha) &= \lam^{k-1}_2(\alpha),
\end{aligned}
\right.
\\
\left\{
\begin{aligned}
\begin{pmatrix}
\dot \y^k_2\\
\dot \lam^k_2
\end{pmatrix}
+ 
\begin{pmatrix}
A & -\nu^{-1} I\\
-I & -A^T
\end{pmatrix}
\begin{pmatrix}
\y^k_2\\
\lam^k_2
\end{pmatrix}
&=
\begin{pmatrix}
0\\
-\hat \y
\end{pmatrix}
\quad \text{ in } (\alpha, T), \\
\y^k_2(\alpha) &= \y^k_1(\alpha), \\
\lam^k_2(T) + \gamma \y^k_2(T) &= \gamma \hat \y(T).
\end{aligned}
\right.
\end{aligned}
\end{equation}
Here, $\y^k_j$ and $\lam^k_j$ represent the restriction of $\y^k$ and
$\lam^k$ to the time subdomain $I_j$, $j = 1, 2$. The parallel version
of this natural Schwarz algorithm~\eqref{eq:SD1} coincides with the
optimized Schwarz algorithm (3a)-(3b) in~\cite{Gander2016}, under the
conditions $p=q=0$ and $\alpha=\beta$ there.

Although algorithm~\eqref{eq:SD1} preserves the forward-backward
structure of the original system~\eqref{eq:forbacksys}, studies
in~\cite{Gander20241, Gander20242} have shown that this structure is
less important for the convergence behavior of Dirichlet--Neumann
and Neumann--Neumann algorithms with time domain
decomposition. Moreover, the forward-backward structure can always be
recovered by using the linear system in~\eqref{eq:forbacksys}, that is
\begin{equation}\label{eq:2id}
\lam = \nu(\dot \y + A \y), \quad \y = \dot \lam - A^T \lam + \hat \y. 
\end{equation}
The two identities~\eqref{eq:2id} also transform a Dirichlet
transmission condition for one state into a particular Robin type
transmission condition for the other state. We can therefore identify
four variants of the classical Schwarz method applied
to~\eqref{eq:forbacksys} with time domain decomposition, as summarized
in Table~\ref{tab:SD} (left).
\begin{table}[t]
\caption{Four non-overlapping Schwarz variants in time. Left:
  Dirichlet transmission conditions. Right: Neumann transmission conditions}
\label{tab:SD}
\centering
\begin{tabular}{ c|c|c|c|c|c|c|c|c|c }
name &SD$_1$ & SD$_2$ & SD$_3$ & SD$_4$ & &SN$_1$ & SN$_2$ & SN$_3$ & SN$_4$  \\
\hline  
$I_1=(0, \alpha)$ & $\lam$ & $\y$ & $\y$ & $\lam$&& $\dot \lam$ & $\dot \y$ & $\dot \y$ & $\dot \lam$ \\  
\hline 
$I_2=(\alpha, T)$ & $\y$ & $\lam$ & $\y$ &$\lam$&& $\dot \y$ & $\dot \lam$ & $\dot \y$ &$\dot \lam$
\end{tabular}

\end{table}

We call these variants SD$_1$ to SD$_4$ (D for Dirichlet) in the first
row. The second row (resp. third row) shows the transmission condition
at the interface of $I_1$ (resp. $I_2$). All four variants use
Dirichlet transmission conditions at the interface, but we can
use~\eqref{eq:2id} to recover the forward-backward structure for
SD$_2$, SD$_3$, and SD$_4$ as explained above.

We now analyze the convergence of these four variants. For simplicity,
we assume that $A$ is symmetric, i.e., $A=A^T\in\R^{n\times n}$. This
allows us to apply a diagonalization, which leads to $n$ independent
$2\times 2$ reduced systems of ODEs. For the algorithm SD$_1$, this
transforms~\eqref{eq:SD1} to,
\begin{equation}\label{eq:SD1reduced}
\begin{aligned}
\left\{
\begin{aligned}
\begin{pmatrix}
\dot z^k_{1, i}\\
\dot \mu^k_{1, i}
\end{pmatrix}
+ 
\begin{pmatrix}
d_i & -\nu^{-1} \\
-1 & -d_i
\end{pmatrix}
\begin{pmatrix}
z^k_{1, i}\\
\mu^k_{1, i}
\end{pmatrix}
&=
\begin{pmatrix}
0\\
-\hat z_i
\end{pmatrix}
\quad \text{ in } (0, \alpha), \\
z^k_{1, i}(0) &= z_{0, i}, \\
\mu^k_{1, i}(\alpha) &= \mu^{k-1}_{2, i}(\alpha),
\end{aligned}
\right.
\\
\left\{
\begin{aligned}
\begin{pmatrix}
\dot z^k_{2, i}\\
\dot \mu^k_{2, i}
\end{pmatrix}
+ 
\begin{pmatrix}
d_i & -\nu^{-1} \\
-1 & -d_i
\end{pmatrix}
\begin{pmatrix}
z^k_{2, i}\\
\mu^k_{2, i}
\end{pmatrix}
&=
\begin{pmatrix}
0\\
-\hat z_i
\end{pmatrix}
\quad \text{ in } (\alpha, T), \\
z^k_{2, i}(\alpha) &= z^k_{1, i}(\alpha), \\
\mu^k_{2, i}(T) + \gamma z^k_{2, i}(T) &= \gamma \hat z_i(T),
\end{aligned}
\right.
\end{aligned}
\end{equation}
where $\bm{z}^k_j := P^{-1}\y^k_j$, $\hat{\bm{z}} := P^{-1}\hat \y$,
$\bm{\mu}^k_j := P^{-1}\lam^k_j$ and $A = PDP^{-1}$ with
$D:=\text{diag}(d_1, \ldots, d_n)$ the eigenvalues of
$A$. Furthermore, $z^k_{j, i}, \hat z_{i}$, and $\mu^k_{j, i}$ denote
the $i$th components of the vectors $\bm{z}^k_j, \hat{\bm{z}}$, and
$\bm{\mu}^k_j$. Note that the assumption of $A$ being symmetric is
only a theoretical tool for the convergence analysis and is not
required to run these algorithms in practice.

To analyze the convergence behavior of the algorithm SD$_1$, we solve
analytically~\eqref{eq:SD1reduced} by eliminating one variable to
obtain a second-order ODE. If we choose to eliminate the adjoint state
$\mu^k_{j, i}$ and use the first identity in~\eqref{eq:2id}, the
Dirichlet transmission condition: $\mu^k_{1, i}(\alpha) =
\mu^{k-1}_{2, i}(\alpha)$ will be transformed into a Robin
transmission condition: $\nu(\dot z^k_{1, i} + d_i z^k_{1,
  i})=\nu(\dot z^{k-1}_{1, i} + d_i z^{k-1}_{1, i})$. Note that
although the natural classical Schwarz algorithm~\eqref{eq:SD1} uses
Dirichlet transmission conditions at the interface, the convergence
analysis actually evaluates a Robin--Dirichlet type algorithm. This
transformation has also been observed in~\cite{Gander2016}. Solving
the resulting second-order ODE allows us to determine the convergence
factor for SD$_1$ as
\begin{equation}\label{eq:rhoSD1}
\rho_{\text{SD}_1}:=\max_{d_i\in D}\Big|\frac{1+\gamma(\sigma_i\coth(b_i)-d_i)}{\nu(\sigma_i\coth(a_i)+d_i)(\sigma_i\coth(b_i)+d_i + \gamma\nu^{-1})}\Big|,
\end{equation}
where $\sigma_i := \sqrt{d_i^2 + \nu^{-1}}$, $a_i = \sigma_i\alpha$ and $b_i = \sigma_i(T-\alpha)$. 
\begin{remark}
Alternatively, one can eliminate the state $z^k_{j, i}$ using the
second identity in~\eqref{eq:2id}, which results in solving a
second-order ODE for $\mu^k_{j, i}$. This approach leads to a
Dirichlet--Robin type algorithm instead, but with the same convergence
factor~\eqref{eq:rhoSD1}, as observed also for Dirichlet--Neumann time
decomposition methods, see~\cite[Appendix A]{Gander20241}.
\end{remark}
To better understand the convergence of the algorithm SD$_1$, we now
study the convergence factor~\eqref{eq:rhoSD1} in detail. We can first
remove the absolute value, since the denominator is positive and
$\sigma_i\coth(b_i)>\sigma_i>d_i$. Next, for a given eigenvalue $d_i$,
we have
\begin{equation}
\label{eq:diff}
\begin{aligned}
&1+\gamma(\sigma_i\coth(b_i)-d_i) - \nu(\sigma_i\coth(a_i)+d_i)(\sigma_i\coth(b_i)+d_i + \gamma\nu^{-1})\\
= &  - \nu d_i^2(\coth(a_i)\coth(b_i) + 1) - (\coth(a_i)\coth(b_i) - 1) - 2\gamma d_i \\
& - \nu\sigma_i d_i(\coth(a_i) + \coth(b_i)) + \gamma\sigma_i(\coth(b_i) - \coth(a_i)).
\end{aligned}
\end{equation}
If $d_i\geq 0$ and $a_i \leq b_i$, then the latter expression is negative, which implies $\rho_{\text{SD}_1}<1$. Hence, we obtain the following result.

\begin{theorem}\label{thm:SD1}
Assume that $A$ is symmetric positive semi-definite (i.e.,
$d_i\geq0$). Then the Schwarz algorithm~\eqref{eq:SD1} converges for
all initial guesses if (i) $\alpha \leq \frac T2$, or (ii) $\gamma=0$.
\end{theorem}
The assumption on the matrix $A$ is natural, for instance, if $A$ is
the finite-difference discretization of the Laplace operator
$-\Delta$. Additionally, setting $\gamma=0$ implies that we are not
considering the final target in~\eqref{eq:model}. In this case, the
convergence factor reads
$\frac{1}{\nu(\sigma_i\coth(a_i)+d_i)(\sigma_i\coth(b_i)+d_i)}$. Taking
the derivative with respect to $d_i$, we find
\[\begin{aligned}
-&\frac{\sigma_i(\coth(a_i)+\coth(b_i)) +2d_i}{\nu(\sigma_i\coth(a_i)+d_i)^2(\sigma_i\coth(b_i)+d_i)^2}-\frac{d_i(\cosh(a_i)\sinh(a_i) -a_i)}{\nu\sigma_i\sinh^2(a_i)(\sigma_i\coth(a_i)+d_i)^2(\sigma_i\coth(b_i)+d_i)}\\
&-\frac{d_i(\cosh(b_i)\sinh(b_i) -b_i)}{\nu\sigma_i\sinh^2(b_i)(\sigma_i\coth(a_i)+d_i)(\sigma_i\coth(b_i)+d_i)^2}.
  \end{aligned}\]
This
derivative is negative if $d_i \geq 0$, since $\cosh(x)\sinh(x) \geq
x$, $\forall x\in \R$. Therefore, we can bound the convergence factor
and find the following result.

\begin{theorem}\label{thm:estimateSD1}
If $A$ is symmetric positive semi-definite and $\gamma=0$, we obtain
the estimate 
\[\rho_{\text{SD}_1} \leq \frac{1}{\nu(\sigma_{\min}
  \coth(\sigma_{\min}\alpha)+d_{\min})(\sigma_{\min}\coth(\sigma_{\min}(T-\alpha))+d_{\min})}
< \frac{1}{\nu(\sigma_{\min}+d_{\min})^2},\] 
where $d_{\min}$ denotes
the smallest eigenvalue of $A$ and $\sigma_{\min} = \sqrt{d_{\min}^2 +
  \nu^{-1}}$.
\end{theorem}
Studying the general form of~\eqref{eq:rhoSD1}, we observe that, for
large eigenvalues $d_i$, the convergence factor approximates
\[\frac{1+\gamma(\sigma_i\coth(b_i)-d_i)}{\nu(\sigma_i\coth(a_i)+d_i)(\sigma_i\coth(b_i)+d_i
  + \gamma\nu^{-1})}\sim_{\infty}\frac1{4\nu d_i^2},\] 
  implying that
high-frequency components converge very fast. For $d_i=0$, the
convergence factor becomes 
\[\rho_{\text{SD}_1}|_{d_i=0} =
\tanh(\sqrt{\nu^{-1}}\alpha)\frac{\gamma\sqrt{\nu^{-1}}\coth(\sqrt{\nu^{-1}}(T-\alpha))+1}{\coth(\sqrt{\nu^{-1}}(T-\alpha))+\gamma\sqrt{\nu^{-1}}},\]
which is close to 1, especially when the control penalization
parameter $\nu$ is small or $\gamma=0$. Hence, low-frequency
components converge very slowly. Based on the monotonicity of
$\rho_{\text{SD}_1}$ with respect to $d_i$ when $\gamma=0$, we can
improve the convergence by introducing a relaxation parameter $\theta$
in the transmission condition to balance the convergence rates of low
and high frequencies. For instance, we can replace the transmission
condition on $I_1$ in~\eqref{eq:SD1} by $\lam^k_1(\alpha) =
f^{k-1}_{\alpha}$ with $ f^{k}_{\alpha} := (1-\theta) f^{k-1}_{\alpha}
+ \theta\lam^{k-1}_2(\alpha)$, and $\theta\in(0, 1)$. The resulting
convergence factor is 
\[\rho^{\theta}_{\text{SD}_1} :=\max_{d_i\in D}|1
-
\theta(1+\frac{1+\gamma(\sigma_i\coth(b_i)-d_i)}{\nu(\sigma_i\coth(a_i)+d_i)(\sigma_i\coth(b_i)+d_i
  + \gamma\nu^{-1})})|.\] 
  Equioscillating between small and large
eigenvalues, we determine the optimal relaxation parameter
\[\theta^*_{\text{SD}_1} :=
\frac{2}{2+\tanh(\sqrt{\nu^{-1}}\alpha)\frac{\gamma\sqrt{\nu^{-1}}\coth(\sqrt{\nu^{-1}}(T-\alpha))+1}{\coth(\sqrt{\nu^{-1}}(T-\alpha))+\gamma\sqrt{\nu^{-1}}}}.\] 
When
$\gamma=0$, the optimal relaxation parameter simplifies to
\[\frac{2}{2+\tanh(\sqrt{\nu^{-1}}\alpha)\tanh(\sqrt{\nu^{-1}}(T-\alpha))},\]
which is approximately $\frac23$ .

For the algorithm SD$_2$, we reverse the transmission conditions
$\mu^k_{1, i}(\alpha) = \mu^{k-1}_{2, i}(\alpha)$ and $z^k_{2,
  i}(\alpha) = z^{k}_{1, i}(\alpha)$ in~\eqref{eq:SD1}, hence also
in~\eqref{eq:SD1reduced}. We thus obtain for SD$_2$ the convergence
factor 
\[\rho_{\text{SD}_2}:=\max_{d_i\in
  D}|\frac{\nu(\sigma_i\coth(a_i)+d_i)(\sigma_i\coth(b_i)+d_i +
  \gamma\nu^{-1})}{1+\gamma(\sigma_i\coth(b_i)-d_i)}|,\] 
  which is the
inverse of $\rho_{\text{SD}_1}$. Hence algorithm SD$_2$ diverges under
the assumption of Theorem~\ref{thm:SD1}, and in particular, it diverges
violently for high-frequency components, because
\[\frac{\nu(\sigma_i\coth(a_i)+d_i)(\sigma_i\coth(b_i)+d_i +
  \gamma\nu^{-1})}{1+\gamma(\sigma_i\coth(b_i)-d_i)}\sim_{\infty}
4\nu d_i^2.\] 
For low-frequency components, the converge is poor,
especially when $\nu$ is small, or it can even diverge when $\gamma=0$, as
\[\rho_{\text{SD}_2}|_{d_i=0} =
\coth(\sqrt{\nu^{-1}}\alpha)\frac{\coth(\sqrt{\nu^{-1}}(T-\alpha))+\gamma\sqrt{\nu^{-1}}}{\gamma\sqrt{\nu^{-1}}\coth(\sqrt{\nu^{-1}}(T-\alpha))+1}.\] 
Based
on these observations, algorithm SD$_2$ is not an efficient algorithm
and can also not be improved with relaxation techniques.

We now study Algorithms SD$_3$ and SD$_4$. Since they pass Dirichlet
data at the interface using only one state, they have similar
behavior, and we just present the analysis for SD$_3$. We replace the
transmission condition $\mu^k_{1, j}(\alpha) = \mu^{k-1}_{2,
  j}(\alpha)$ on $I_1$ with $z^k_{1, j}(\alpha) = z^{k-1}_{2,
  j}(\alpha)$ in~\eqref{eq:SD1reduced}.  Using the second-order ODE
for $z^k_{j, i}$, we find that the transmission conditions still
remain, i.e., $z^k_{1, j}(\alpha) = z^{k-1}_{2, j}(\alpha)$ on $I_1$
and $z^k_{2, j}(\alpha) = z^{k}_{1, j}(\alpha)$ on $I_2$. This
indicates that this is a Dirichlet--Dirichlet type algorithm in the
classical Schwarz sense, thus suffering from the same non-convergence
as the classical Schwarz algorithm without overlap. Indeed, its
convergence factor $\rho_{\text{SD}_3}$ equals 1 for all eigenvalues
$d_i$, and this cannot be improved with relaxation. Similarly, using
the second-order ODE for $\mu^k_{j, i}$, we obtain also
$\rho_{\text{SD}_4}=1$. Hence, in contrast to the Dirichlet--Neumann
algorithms in time \cite{Gander20241}, among all four Schwarz variants
with Dirichlet transmission conditions, only algorithm
SD$_1$~\eqref{eq:SD1}, which naturally preserves the forward-backward
structure, exhibits good convergence behavior.

\section{Neumann transmission condition}\label{sec:3}

Schwarz methods with Neumann transmission conditions are not used for
elliptic problems, since they are not convergent in general, as one
can see from a simple 1D example. We investigate now if Schwarz
methods in time for parabolic optimal control problems can be promising solvers.
We identity once again four variants, see Table~\ref{tab:SD} (right),
called SN$_1$ to SN$_4$. Similar to SD$_1$, algorithm SN$_1$ naturally
retains the forward-backward structure. For the other three variants,
this structure can be recovered using identities in~\eqref{eq:2id}.

To analyze the convergence behavior of the four variants, we follow
the same approach as in Section~\ref{sec:2}. Algorithm SN$_1$ is
similar to~\eqref{eq:SD1}, but with transmission conditions replaced
by $\dot\lam^k_{1, i}(\alpha) = \dot\lam^{k-1}_{2, i}(\alpha)$ on
$I_1$ and $\dot \y^k_{2, i}(\alpha) = \dot \y^{k}_{1, i}(\alpha)$ on
$I_2$. When analyzing its convergence using the second-order ODE for
$z^k_{i, j}$, we are then examining a Robin--Neumann type algorithm.
We find for SN$_1$ the convergence factor
\begin{equation}\label{eq:SN1}
\rho_{\text{SN}_1}:=\max_{d_i\in D}\Big|\frac{1 + \gamma(\sigma_i\tanh(b_i) - d_i)}{\nu(\sigma_i\tanh(a_i) + d_i)(\sigma_i\tanh(b_i) + d_i + \gamma\nu^{-1})}\Big|,
\end{equation}
similar to $\rho_{\text{SD}_1}$, with hyperbolic cotangent functions
in~\eqref{eq:rhoSD1} replaced by hyperbolic tangent
functions. However, unlike Theorem~\ref{thm:SD1}, we cannot directly
obtain a similar result for SN$_1$, since the sign of~\eqref{eq:diff}
is less clear when replacing hyperbolic cotangent by hyperbolic
tangent. Nevertheless, substituting $\gamma=0$ into~\eqref{eq:SN1}
yields $\frac{1}{\nu(\sigma_i\tanh(a_i) + d_i)(\sigma_i\tanh(b_i) +
  d_i)}$, which is a decreasing function of $d_i$, as
$\sigma_i$ and the hyperbolic tangent are both increasing functions of
$d_i$ when $d_i\geq 0$. Therefore, we obtain a similar result as
Theorem~\ref{thm:estimateSD1}.

\begin{theorem}\label{thm:estimateSN1}
If $A$ is symmetric positive semi-definite and $\gamma=0$, we obtain the estimate 
\[\rho_{\text{SN}_1} \leq  \frac{1}{\nu(\sigma_{\min} \tanh(\sigma_{\min}\alpha)+d_{\min})(\sigma_{\min}\tanh(\sigma_{\min}(T-\alpha))+d_{\min})}.\]
\end{theorem}
Moreover, for large eigenvalues, the convergence factor for SN$_1$ is
approximately 
\[\frac{1 + \gamma(\sigma_i\tanh(b_i) -
  d_i)}{\nu(\sigma_i\tanh(a_i) + d_i)(\sigma_i\tanh(b_i) + d_i +
  \gamma\nu^{-1})}\sim_{\infty}\frac1{4\nu d_i^2},\] 
  meaning that
SN$_1$ is a very good smoother for high-frequency components. For a
zero eigenvalue $d_i=0$, 
\[\rho_{\text{SN}_1}|_{d_i=0} =
\coth(\sqrt{\nu^{-1}}\alpha)\frac{\gamma\sqrt{\nu^{-1}}\tanh(\sqrt{\nu^{-1}}(T-\alpha))+1}{\tanh(\sqrt{\nu^{-1}}(T-\alpha))+\gamma\sqrt{\nu^{-1}}}.\] 
Thus,
low-frequency components converge very slowly, especially when $\nu$
is small, or it can diverge when $\gamma=0$.
As with SD$_1$, one can use the monotonicity of
$\rho_{\text{SN}_1}$ in the case $\gamma=0$ and improve the
convergence with a relaxation parameter $\theta$. The convergence
factor with relaxation is 
\[\rho^{\theta}_{\text{SN}_1}:=\max_{d_i\in
  D} |1-\theta(\frac{1 + \gamma(\sigma_i\tanh(b_i) -
  d_i)}{\nu(\sigma_i\tanh(a_i) + d_i)(\sigma_i\tanh(b_i) + d_i +
  \gamma\nu^{-1})})|.\] 
  Using the equioscillation principle, we
determine the optimal relaxation parameter 
\[\theta^*_{\text{SN}_1} :=
\frac{2}{2+\coth(\sqrt{\nu^{-1}}\alpha)\frac{\gamma\sqrt{\nu^{-1}}\tanh(\sqrt{\nu^{-1}}(T-\alpha))+1}{\tanh(\sqrt{\nu^{-1}}(T-\alpha))+\gamma\sqrt{\nu^{-1}}}}.\] 
For
$\gamma=0$, this becomes
$\frac{2}{2+\coth(\sqrt{\nu^{-1}}\alpha)\coth(\sqrt{\nu^{-1}}(T-\alpha))}$
and is also bounded by $\frac23$.

For algorithm SN$_2$, we reverse the transmission condition in SN$_1$
and obtain again the inverse of the convergence factor of SN$_1$
in~\eqref{eq:SN1}. As for SD$_2$, SN$_2$ is hence not an efficient
algorithm and cannot be improved using relaxation. For algorithms
SN$_3$ and SN$_4$, they pass Neumann data at the interface using only
one state. Similarly as for SD$_3$ and SD$_4$, we find that
$\rho_{\text{SN}_3} = \rho_{\text{SN}_4} = 1$ for all eigenvalues
$d_i$, indicating stagnation and no improvement with
relaxation. Hence, among the four variants with Neumann transmission
conditions, only algorithm SN$_1$, which naturally preserves the
forward-backward structure, has good convergence behavior, and this
even though it is a Schwarz method with Neumann transmission
conditions, which do not work in the elliptic case!

\section{Numerical experiments and comments}\label{sec:4}

\begin{figure}
\label{fig:test}
\centering
\includegraphics[scale=0.13]{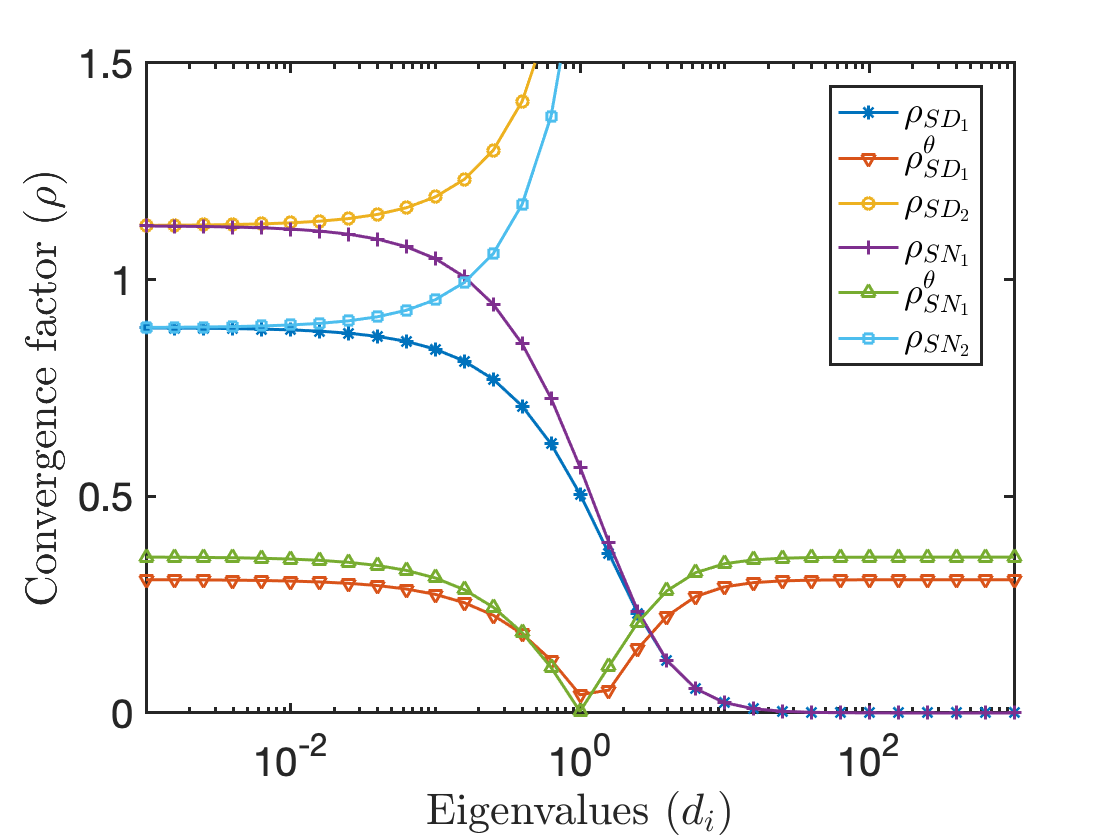}
\includegraphics[scale=0.13]{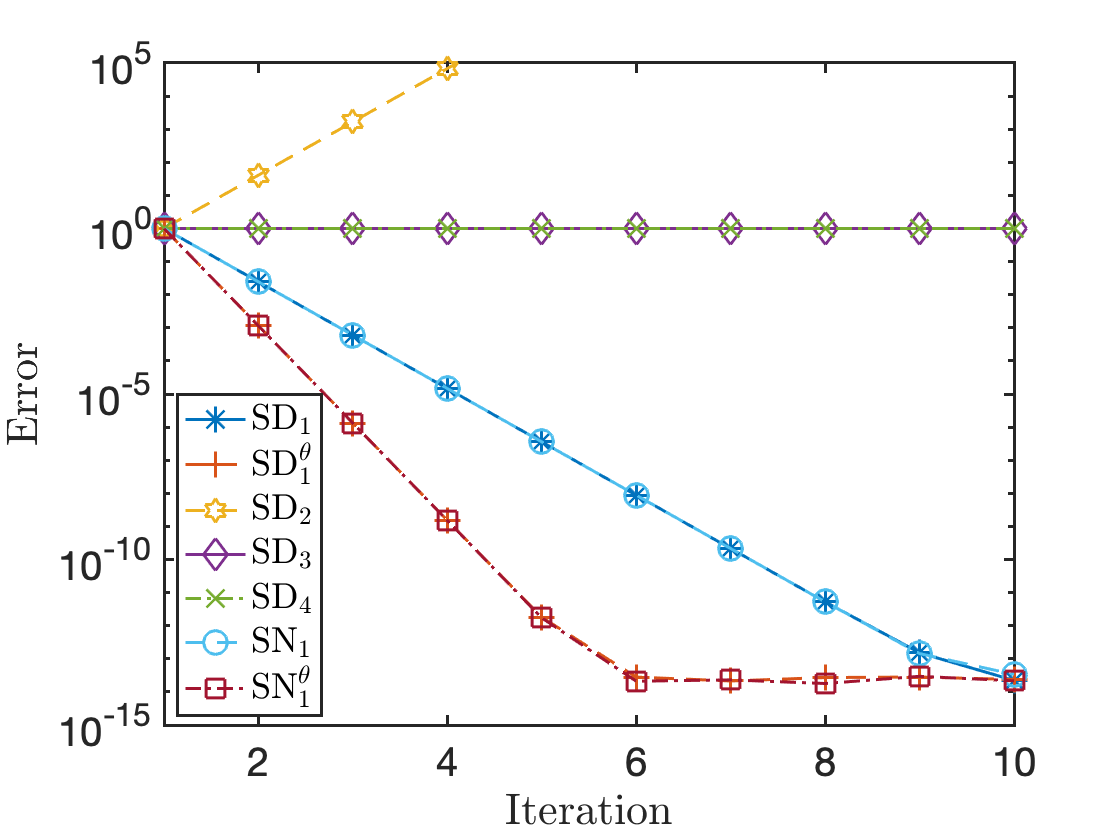}
\caption{Convergence factor as a function of eigenvalues (left) and error decay as a function of the number of iterations (right).}
\end{figure}

We first plot the convergence factor $\rho$ as a function of the eigenvalues $d_i$ in Figure~\ref{fig:test} (left). 
We set the parameters $\nu=0.1$, $\gamma=10$, $T=1$ and $\alpha=0.4$, and observe that both SD$_2$ and SN$_2$ diverge for $d_i \geq1$. Algorithms SD$_1$ and SN$_1$ are two good smoothers for high-frequency components, but they exhibit poor convergence for low-frequency components. This is significantly improved when using relaxation techniques. We find numerically $\theta^*_{\text{SD}_1}\approx0.692$ and $\theta^*_{\text{SN}_1}\approx0.640$, which are consistent with their theoretical values. To evaluate the numerical performance of these variants, we apply them to solve the one-dimensional heat control problem $\partial_t y - \partial_{xx} y = u$ with homogenous Dirichlet boundary conditions and a zero initial condition. We keep the same parameter values and choose the target state $\hat y=\sin(\pi x)(2t^2 + t)$. We use the Crank-Nicolson scheme with mesh size $h_t = h_x = 1/32$. The error decay as a function of the number of iterations is shown in Figure~\ref{fig:test} (right). As expected, SD$_2$ diverge violently, and both SD$_3$ and SD$_4$ stagnate. The convergence of SD$_1$ and SN$_1$ is already efficient without relaxation, due to the smallest eigenvalue in this test case being around 10. The convergence can be improved from 10 to 6 iterations with a relaxation parameter $\theta=0.975$ for both algorithms.

Unlike our observation in~\cite{Gander20241, Gander20242} for
Dirichlet--Neumann and Neumann--Neumann algorithms in time, classical
Schwarz algorithms with only Dirichlet or Neumann transmission
conditions are much more sensitive to the forward-backward
structure. We observe that only SD$_1$ and SN$_1$, which naturally
preserve this structure, have good convergence behavior. All other
variants preform poorly and cannot be improved even with relaxation
techniques. For SD$_1$ and SN$_1$, we also provided estimates and
closed-form expressions for the optimal relaxation
parameters. In~\cite{Gander2016}, the authors used transmission
conditions of the form $\lam+p\y$ on $I_1$ and $\y - q \lam$ on $I_2$,
with two parameters $p, q\geq 0$. It would be interesting to extend
this approach using transmission conditions of the form $\dot
\lam+p\dot \y$ and $\dot \y - q \dot \lam$ to further improve the
convergence in the case of Neumann transmission conditions only.

\bibliography{biblio}

\begin{thebibliography}{1}

\bibitem{Gander2008}
M.~J. Gander.
\newblock Schwarz methods over the course of time.
\newblock {\em Electronic Transactions on Numerical Analysis}, 31:228--255,
  2008.

\bibitem{Gander2016}
M.~J. Gander and F.~Kwok.
\newblock Schwarz methods for the time-parallel solution of parabolic control
  problems.
\newblock In T.~Dickopf, M.~J. Gander, L.~Halpern, R.~Krause, and L.~F.
  Pavarino, editors, {\em Domain Decomposition Methods in Science and
  Engineering XXII}, pages 207--216, Cham, 2016. Springer.

\bibitem{Gander20241}
M.~J. Gander and L.-D. Lu.
\newblock New time domain decomposition methods for parabolic optimal control
  problems {I}: {D}irichlet--{N}eumann and {N}eumann--{D}irichlet algorithms.
\newblock To appear in SIAM Journal on Numerical Analysis, 2024.

\bibitem{Gander20242}
M.~J. Gander and L.-D. Lu.
\newblock New time domain decomposition methods for parabolic optimal control
  problems {II}: {N}eumann--{N}eumann algorithms.
\newblock Submitted, 2024.

\bibitem{Lions1990}
P.-L. Lions.
\newblock On the {S}chwarz {A}lternating {M}ethod {III}: {A} {V}ariant for
  {N}onoverlapping {S}ubdomains.
\newblock In T.~F. Chan, R.~Glowinski, J.~Périaux, and O.~Widlund., editors,
  {\em Third International Symposium on Domain Decomposition Methods for
  Partial Differential Equations}, pages 202--223. SIAM, 1990.

\bibitem{Schwarz1870}
H.~A. Schwarz.
\newblock Über einen {G}renzübergang durch alternierendes {V}erfahren.
\newblock {\em Vierteljahrsschrift der {N}aturforschenden {G}esellschaft in
  {Z}ürich}, 15:272--286, 1870.

\end{thebibliography}
\bibliographystyle{plain}
\end{document}